\input amstex
\documentstyle{amsppt}
\pagewidth{5.9in}\vsize8.5in\parindent=6mm
\parskip=4pt\baselineskip=6pt\tolerance=10000\hbadness=500
\magnification=1100

\rightheadtext{Smooth cone type multipliers}
\document
\topmatter
\title Weak type estimates on certain Hardy spaces \\
 for smooth cone type multipliers
\endtitle
\author Sunggeum Hong and Yong-Cheol Kim  \endauthor
\abstract Let $\varrho\in C^{\infty} ({\Bbb R}^d\setminus\{0\})$
be a non-radial homogeneous distance function satisfying
$\varrho(t\xi)=t\varrho(\xi)$. For $f\in\frak S ({\Bbb R}^{d+1})$
and $\delta>0$, we consider convolution operator ${\Cal
T}^{\delta}$ associated with the smooth cone type multipliers
defined by
$$\widehat {{\Cal T}^{\delta} f}(\xi,\tau)=
\left(1-\frac{\varrho(\xi)}{|\tau|} \right)^{\delta}_+\hat f
(\xi,\tau),\,\,(\xi,\tau)\in {\Bbb R}^d \times \Bbb R.$$ If the
unit sphere $\Sigma_{\varrho}\fallingdotseq\{\xi\in {\Bbb R}^d :
\varrho(\xi)=1\}$ is a convex hypersurface of finite type, then we
prove that the operator ${\Cal T}^{\delta(p)}$ maps from
$H^p({\Bbb R}^{d+1})$, $0<p<1$, into weak-$L^p({\Bbb R}^{d+1})$
for the critical index $\delta(p)=d(1/p -1/2)-1/2$.
\endabstract
\thanks 2000 Mathematics Subject Classification: 42B15, 42B30.
\endthanks
\thanks The second author was supported in part by a Korea University Grant.
\endthanks
\address Department of Mathematics, Chosun University, Gwangju 501-759, Korea
\endaddress
\email skhong$\@$mail.chosun.ac.kr \endemail
\address Department of Mathematics Ed., Korea University, Seoul 136-701, Korea
\endaddress
\email ychkim$\@$korea.ac.kr \endemail
\endtopmatter

\define\supp{{\text{\rm supp }}}

\define\inn#1#2{\langle#1,#2\rangle}

\define\lcontr{\rfloor}
\define\lco#1#2{{#1}\lcontr{#2}}
\define\lcoi#1#2{\imath({#1}){#2}}
\define\rco#1#2{{#1}\rcontr{#2}}
\redefine\exp{{\text{\rm exp}}}

\define\ap{\alpha}             

\define\gm{\gamma}             \define\Gm{\Gamma}
\define\dt{\delta}             
\define\vep{\varepsilon}
\define\zt{\zeta}
\define\th{\theta}

\define\ld{\lambda}            
\define\sm{\sigma}             \define\Sm{\Sigma}
\define\vp{\varphi}
             \define\Om{\Omega}
\define\vr{\varrho}            \define\iy{\infty}
\define\lt{\left}            \define\rt{\right}
\define\f{\frac}             \define\el{\ell}

\define\fJ{{\frak J}}

\define\fS{{\frak S}}

\define\fa{{\frak a}}
\define\fb{{\frak b}}

\define\fe{{\frak e}}

\define\fg{{\frak g}}
\define\fh{{\frak h}}

\define\fk{{\frak k}}

\define\fn{{\frak n}}

\define\BN{{\Bbb N}}

\define\BR{{\Bbb R}}

\define\BZ{{\Bbb Z}}

\define\cA{{\Cal A}}
\define\cB{{\Cal B}}
\define\cC{{\Cal C}}
\define\cD{{\Cal D}}
\define\cE{{\Cal E}}
\define\cF{{\Cal F}}

\define\cI{{\Cal I}}
\define\cJ{{\Cal J}}
\define\cK{{\Cal K}}

\define\cO{{\Cal O}}
\define\cP{{\Cal P}}
\define\cQ{{\Cal Q}}
\define\cR{{\Cal R}}

\define\cT{{\Cal T}}
\define\cU{{\Cal U}}
\define\cV{{\Cal V}}




\define\la{\langle}          \define\ra{\rangle}
\define\s{\setminus}         
\define\n{\nabla}            
\define\pa{\partial}        \define\fd{\fallingdotseq}
       
     \define\ds{\dsize}
        \define\ls{\lesssim }
\define\gs{\gtrsim}

\subheading{1. Introduction}

Let $M$ be a real-valued $d\times d$ matrix whose eigenvalues have
positive real part. Then the linear transformations
$A_t=\exp(M\,\log t)$, $t>0$, form a dilation group on $\BR^d$.
Let $\varrho\in C^{\iy}(\BR^d\s\{0\})$ be a positive real-valued
function satisfying $\varrho(A_t\,\xi)=t\,\varrho(\xi)$, which is
called an $A_t$-homogeneous distance function ( refer to [11] for
elementary properties ).

Let $\fS(\BR^{d+1})$ be the Schwartz space of rapidly decreasing
$C^{\iy}(\BR^{d+1})$-functions, and let $\hat f$ be the Fourier
transform of $f\in\fS(\BR^{d+1})$. In what follows, we shall
always assume that $A_t=t\,I$ and $\varrho\in
C^{\iy}(\BR^d\s\{0\})$ is a non-radial $A_t$-homogeneous distance
function whose unit sphere $$\Sm_{\vr}\fd\{\xi\in\BR^d :
\vr(\xi)=1\}$$ satisfies a finite type condition, i.e. every
tangent line makes finite order of contact with $\Sm_{\vr}$. For
$f\in\fS (\BR^{d+1})$ and $\dt>0$, we define convolution operators
$\cT^{\dt}$ by
$$\cT^{\dt} f(x,t)=\f{1}{(2\pi)^{d+1}}\int_{\BR}\int_{\BR^d}e^{i\la
x,\xi \ra +it\tau} \lt(1-\f{\vr(\xi)}{|\tau|}\rt)^{\dt}_+ \hat f
(\xi,\tau) \,d\xi\,d\tau,\,\,(\xi,\tau)\in\BR^d\times\BR.$$

In the case that $A_t=t^{1/2}\,I$ and $\varrho(\xi)=|\xi|^2$,
there is no optimal $L^p(\BR^d)$-mapping properties of $\cT^{\dt}$
which are known for $p>1$ and $d\ge 2$. For partial results
related with this, the reader can refer to G. Mockenhaupt [7] and
J. Bourgain [1]. However, it is known ( see [5] ) about its sharp
weak type results on $H^p(\BR^{d+1})$, $0<p<1$, with the critical
index $\dt(p)=d(1/p-1/2)-1/2$.

The purpose of this article is to obtain sharp weak type endpoint
($\dt(p)=d(1/p-1/2)-1/2$) results of $\cT^{\dt}$ on
$H^p(\BR^{d+1})$, $0<p<1$, for the distance function $\vr(\xi)$ in
the sense that $\cT^{\dt}$ is a bounded operator of
$H^p(\BR^{d+1})$ into $L^p(\BR^{d+1})$ for $\dt>\dt(p)$, but it
fails to be of restricted weak type $(p,p)$ on $H^p(\BR^{d+1})$
for all $\dt<\dt(p)$ and fails to be a bounded operator of
$H^p(\BR^{d+1})$ into $L^p(\BR^{d+1})$ when $\dt=\dt(p)$. Here
$H^p$ denotes the standard real Hardy space as defined by E. M.
Stein in [9].

\proclaim{Theorem 1.1} Suppose that $\vr\in C^{\iy}(\BR^d\s\{0\})$
is a non-radial homogeneous distance function satisfying
$\vr(t\xi)=t\vr(\xi)$ whose unit sphere $\Sm_{\vr}$ is a convex
hypersurface of finite type. If $\dt(p)=d(1/p-1/2)-1/2$ for
$0<p<1$, then $\cT^{\dt(p)}$ maps $H^p(\BR^{d+1})$ boundedly into
weak-$L^p(\BR^{d+1})$; that is, there exists a constant $C>0$ not
depending upon $\ld$ and $f$ such that for any $f\in
H^p(\BR^{d+1})$,
$$\lt|\{(x,t)\in\BR^{d+1} : |\cT^{\dt(p)} f(x,t)|>\ld\}\rt|\le
C\, \f{\|f\|^p_{H^p(\BR^{d+1})}}{\ld^p},\,\,\ld>0.$$
\endproclaim

\noindent {\it Remarks.} $(i)$ This result generalizes that of [5]
in $\BR^d\times\BR$ to non-radial cases ( which is of finite type
and convex ) by using the arguments based on the results in [6].

$(ii)$ As a matter of fact, we prove this result under more
general surface condition than the finite type condition on
$\Sm_{\vr}$, which was so called a spherically integrable
condition of order $<1$ in [6].

In what follows, we shall use the polar coordinates; given
$x\in\BR^d$, we write $x=r\,\th$ where $r=|x|$ and
$\th=(\th_1,\th_2,\cdots,\th_d)\in S^{d-1}$. Given two quantities
$A$ and $B$, we write $A\ls B$ or $B\gs A$ if there is a positive
constant $c$ ( possibly depending on the dimension $d\ge 2$, the
hypersurface $\Sm_{\vr}$, and the index $p$ to be given ) such
that $A\le c B$. We also write $A\sim B$ if $A\ls B$ and $B\ls A$.
We denote by $\BR^d_0\fd \BR^d\s\{0\}$. For $e\in S^{d-1}$, we let
$\cD_e f$ denote the directional derivative $\la e,\n_{\xi}\ra f$
of a function $f(\xi)$ defined on $\BR^d$. For $\fe\in S^{d-2}$,
we let $\cD^w_{\fe} g$ denote the directional derivative
$\la\fe,\n_w\ra g$ of a function $g(w)$ defined on $\BR^{d-1}$.
For a multi-index $\ap=(\ap_1,\ap_2,\cdots,\ap_{d+1})\in
(\BN\cup\{0\})^{d+1}$, we denote its size $|\ap|$ by $|\ap|\fd
\ap_1+\ap_2+\cdots +\ap_{d+1}$.

\subheading{2. Preliminary estimates on a smooth convex
hypersurface of finite type}

Let $\Sm$ be a smooth convex hypersurface of $\BR^d$ and let
$d\sm$ be the induced surface area measure on $\Sm$. For
$x\in\BR^d$, we denote by $\cB (\xi(x),s)$ the spherical cap near
$\xi(x)\in\Sm$ cut off from $\Sm$ by a plane parallel to
$T_{\xi(x)}(\Sm)$ ( the affine tangent plane to $\Sm$ at $\xi(x)$
) at distance $s>0$ from it; that is, $$\cB
(\xi(x),s)=\{\xi\in\Sm:\,d(\xi,T_{\xi(x)}(\Sm))<s\},$$ where
$\xi(x)$ is the point of $\Sm$ whose outer unit normal is in the
direction $x$. These spherical caps play an important role in
furnishing the decay of the Fourier transform of the measure
$d\sm$. It is well known [9] that the function
$$\varPhi(\th)\fd\sup_{r>0}\,\sm[\cB(\xi(r\th),1/r)](1+r)^{\f{d-1}{2}}\tag{2.1}$$
is bounded on $S^{d-1}$ provided that $\Sm$ has nonvanishing
Gaussian curvature.

\noindent{\it Remark.} (i) B. Randol [8] proved that if $\Sm$ is a
real analytic convex hypersurface of $\BR^d$ then $\varPhi\in
L^p(S^{d-1})$ for some $p>2$. Thus any real analytic convex
hypersurface satisfies that $\varPhi\in L^p(S^{d-1})$ for any
$p\le 2$.

\noindent\quad(ii) More generally, it was shown by I. Svensson
[12] that if $\Sm$ is a smooth convex hypersurface of finite type
$k\ge 2$ then $\varPhi\in L^p(S^{d-1})$ for some $p>2$.

Thus, by the above remark (ii), it is natural for us to obtain the
following corollary.

\proclaim{Corollary 2.1} If $\Sm$ is a smooth convex hypersurface
of $\BR^d$ which is of finite type, then $\varPhi\in L^p(S^{d-1})$
for any $p\le 2$.
\endproclaim

Sharp decay estimates for the Fourier transform of surface measure
on a smooth convex hypersurface $\Sm$ of finite type $k\ge 2$ have
been obtained by Bruna, Nagel, and Wainger [2]; precisely
speaking,
$$\lt|\cF[d\sm](x)\rt|\sim\sm[\cB(\xi(x),1/|x|)].\tag{2.2}$$ They
define a family of anisotropic balls on $\Sm$ by letting
$$\cB (\xi_0,s)=\{\xi\in\Sm:\,d(\xi,T_{\xi_0}(\Sm))<s\}$$ where
$\xi_0\in \Sm$. We now recall some properties of the anisotropic
balls $\cB(\xi_0,s)$ associated with $\Sm$. The proof of the
doubling property in [2] makes it possible to obtain the following
stronger estimate for the surface measure of these balls; $$
\sm[\cB(\xi_0,\gm \,s)]\ls \lt\{ \aligned \gm^{\f{d-1}{2}}
\sm[\cB(\xi_0,s)],\,\,\,\,&\gm\ge 1, \\ \gm^{\f{d-1}{k}}
\sm[\cB(\xi_0,s)],\,\,\,\,&\gm< 1. \endaligned \rt. $$ It also
follows from the triangle inequality and the doubling property [2]
that there is a positive constant $C>0$ independent of $s>0$ such
that $$\f{1}{C}\,\sm[\cB(\xi_0,s)]\le \sm[\cB(\xi,s)]\le
C\,\sm[\cB(\xi_0,s)]\,\,\,\text { for any $\xi\in\cB(\xi_0,s)$.}
\tag{2.3}$$

\proclaim{Lemma 2.2 [6]} Let $\Sm$ be a smooth convex hypersurface
of $\BR^d$ which is of finite type $k\ge 2$. Then there is a
constant $C=C(\Sm)>0$ such that for any $y\in B(0;s)$ and $x\in
{B(0;2s)}^c$, $0<s\le 1$, $$\xi(x-y)\in\cB(\xi(x),C/|x|)$$ where
$\xi(x)$ is the point of $\Sm$ whose outer unit normal is in the
direction $x$.
\endproclaim

\proclaim{Lemma 2.3} Let $\Sm$ be a smooth convex hypersurface of
$\BR^d$ which is of finite type $k\ge 2$. Then there is a constant
$C=C(\Sm)>0$ such that for any  $x,y\in\BR^d$ with $|x|>2\,|y|>0$,
$$\varPhi\lt(\f{x-y}{|x-y|}\rt)\le
C\,\,\varPhi\lt(\f{x}{|x|}\rt)$$ where $\varPhi$ is the function
defined as in $(2.1)$.
\endproclaim

\noindent{\it Proof.} It easily follows from (2.3), the definition
of $\varPhi$, and Lemma 2.2 that for any $y\in B(0;s)$ and $x\in
{B(0;2s)}^c$, $0<s\le 1$, $$ \split
\varPhi\lt(\f{x-y}{|x-y|}\rt)&\ls\sup_{r>0}\sm[\cB(\xi(x-y),1/r)]\,
(1+r)^{\f{d-1}{2}} \\
&\ls \sup_{r>0}\sm[\cB(\xi(x),1/r)]\,
(1+r)^{\f{d-1}{2}}=\varPhi\lt(\f{x}{|x|}\rt). \endsplit
$$ Thus this implies that for any  $x,y\in\BR^d$ with $|x|>2\,|y|>0$,
$$\varPhi\lt(\f{x-y}{|x-y|}\rt)
=\varPhi\lt(\f{x/|y|-y/|y|}{|x/|y|-y/|y||}\rt)\ls\varPhi\lt(\f{x/|y|}{|x/|y||}\rt)
=\varPhi\lt(\f{x}{|x|}\rt). \,\,\,\qed$$

\subheading{3. The kernel estimate associated with the cone type
multipliers }

Throughout this section from now on, we shall concentrate upon
obtaining decay estimate of the kernel associated with the cone
type multipliers. First of all, we obtain the lower bound of the
phase function of the kernel on the dual cone
$\Gm_{\gm}=\{(x,t)\in\BR^d\times\BR : |t|\ge\gm |x|\},
\gm=\sup_{\xi\in\cB_{\vr}}\,|\xi|$, and the unit ball
$\cB_{\vr}\fd\{\xi\in\BR^d : \vr(\xi)\le 1\}$.

\proclaim{Lemma 3.1} Suppose that $\vr\in C^{\iy}(\BR^d\s\{0\})$
is a non-radial homogeneous distance function satisfying
$\vr(t\xi)=t\vr(\xi)$. Then we have the following estimate
$$\inf_{\xi\in\cB_{\vr}}\, |t+\la x,\xi
\ra|\ge |t|-\gm |x|\ge 0,\,\,\,(x,t)\in\Gm_{\gm}.$$
\endproclaim

\noindent{\it Proof.} It easily follows from the triangle
inequality and Schwarz inequality. \qed

Let $\psi\in C_0^{\iy}(\BR)$ be supported in $(1/2,2)$ such that
$\sum_{l\in\BZ}\psi(2^{-l}\,t)=1$ for $t>0$. For fixed $l\in\BZ$,
we shall need pointwise estimates for the operators
$$\cT_l^{\dt(p)}
f(x,t)=\f{1}{(2\pi)^{d+1}}\int_{\BR}\int_{\BR^d}\cK^{\dt(p)}_l(x-y,t-s)\,f(y,s)\,dy\,ds$$
where
$$\cK^{\dt(p)}_l(x,t)=\f{1}{(2\pi)^{d+1}}\int_{\BR}\int_{\BR^d}e^{i\la x,\xi \ra +it\tau}
\lt(1-\f{\vr(\xi)}{|\tau|}\rt)^{\dt(p)}_+
\psi(2^{-l}\,\tau)\,d\xi\,d\tau. \tag{3.1}$$ For each $l\in\BZ$,
the kernel $\cK^{\dt(p)}_l$ has the property
$$\cK^{\dt(p)}_l(x,t)=2^{(d+1)l}\,\cK^{\dt(p)}_0(2^l x,2^l t).
\tag{3.2}$$ Let $\phi\in C_0^{\iy}(\BR)$ be supported in $(1/2,2)$
such that $\sum_{k\in\BZ}\phi(2^k\,s)=1$ for $0<s<1$, and let
$$ \cK^{\dt(p)}_{k,l}(x,t)=\f{1}{(2\pi)^{d+1}}\int_{\BR}\int_{\BR^d}e^{i\la x,\xi \ra +it\tau}
\phi\lt(2^{k+1}\lt(1-\f{\vr(\xi)}{|\tau|}\rt)\rt)
\lt(1-\f{\vr(\xi)}{|\tau|}\rt)^{\dt(p)}_+
\psi(2^{-l}\,\tau)\,d\xi\,d\tau. $$ We write
$\sum_{k\in\BN}\cK^{\dt(p)}_{k,l}+\cK^{\dt(p)}_{0,l}=\cK^{\dt(p)}_l$
and
$\sum_{k\in\BN}\cT^{\dt(p)}_{k,l}+\cT^{\dt(p)}_{0,l}=\cT^{\dt(p)}_l$.
Since the kernel $\cK^{\dt(p)}_{0,l}(x,t)$ has a very nice decay
in terms of $x$-variable and a properly good decay in
$t$-variable, it suffices to treat only the operators
$\cT^{\dt(p)}_{k,l}$ in the following lemma.

\proclaim{Lemma 3.2} Suppose that $\vr\in C^{\iy}(\BR^d\s\{0\})$
is a non-radial homogeneous distance function satisfying
$\vr(t\xi)=t\vr(\xi)$. If $|x|\le 2^{2-l}\gm^{-1}$ for fixed
$l\in\BZ$, then for each $k\in\BN$ we have the following uniform
estimates; for any $N\in\BN$,
$$|\cK^{\dt(p)}_{k,l}(x,t)|+\sum_{|\ap|=n+1}\f{1}{\ap
!}\,|[\cD^{\ap}\cK^{\dt(p)}_{k,l}](x,t)|\le
C\,2^{(d+1)l}2^{-k(\dt(p)+1)}\min\{1,(2^{-k}2^l\,|t|)^{-N}\},$$
and thus $$|\cK^{\dt(p)}_l(x,t)|+\sum_{|\ap|=n+1}\f{1}{\ap
!}\,|[\cD^{\ap}\cK^{\dt(p)}_l](x,t)|\le
C\,2^{(d+1)l}\f{1}{(1+2^l\,|t|)^{\dt(p)+1}} $$ where $n\in\BN$,
$\ap=(\ap_1,\ap_2,\cdots,\ap_{d+1})\in (\BN\cup\{0\})^{d+1}$ is a
multi-index, and we denote by $\ap !=\ap_1 ! \,\ap_2 !\cdots
\,\ap_{d+1} !$ and
$$\cD^{\ap}\cK^{\dt(p)}_{k,l}(x,t)=\lt(\f{\pa}{\pa x_1}\rt)^{\ap_1}\lt(\f{\pa}{\pa x_2}\rt)^{\ap_2}
\cdots\lt(\f{\pa}{\pa x_d}\rt)^{\ap_d}\lt(\f{\pa}{\pa
t}\rt)^{\ap_{d+1}}\cK_{k,l}(x,t).$$
\endproclaim

\noindent{\it Proof.} It easily follows from the integration by
parts $N$-times with respect to $\tau$. \qed

Now we proceed the case that $|x|>2^{2-l}\gm^{-1}$. By the change
of variables, the integral (3.1) becomes
$$\cK^{\dt(p)}_l(x,t)=2^{(d+1)l}\int_{\BR}\int_{\BR^d} e^{i 2^l|x|\la \f{x}{|x|},\xi\ra\tau+i 2^l t\tau}
(1-\vr(\xi))^{\dt(p)}_+ \,\psi(\tau)\,\tau^2\,d\xi\,d\tau.$$ We
shall employ a decomposition of the Bochner-Riesz multiplier
$(1-\vr)^{\dt(p)}_+$ as in A. C\'ordoba [3].  Let $\phi\in
C^{\iy}_0(1/2,2)$ satisfy $\sum_{k\in\BZ} \phi(2^k t)=1$ for
$t>0$. Put
$\Phi^{\dt(p)}_k=\phi(2^{k+1}(1-\vr))(1-\vr)_+^{\dt(p)}$ and
$\Phi^{\dt(p)}_0=(1-\vr)^{\dt(p)}_+-\sum_{k\in\BN}\Phi^{\dt(p)}_k$
for $k\in\BN$. Then we note that
$\sum_{k\in\BN}\Phi^{\dt(p)}_k=\vp\,(1-\vr)^{\dt(p)}_+\,\,a.e.$,
where $\vp\in C_0^{\iy}(\BR^d)$ is a function supported in the
closed annulus $$\{\xi\in\BR^d :\,1/2\le\vr(\xi)\le 2\}$$ such
that
$$\vp(\xi)=\sum_{k\in\BN}\phi(2^{k+1}(1-\vr(\xi)))\,\,\text{ on the open annulus
$\{\xi\in\BR^d:1/2<\vr(\xi)<1\}$}.$$ We now introduce a partition
of unity $\Xi_{\el}, \,\el=1,2,\cdots, N_0$, on the unit sphere
$\Sm_{\vr}$ which we extend to $\BR^d$ by way of
$\Pi_{\el}(t\,\zt)=\Xi_{\el}(\zt),\,t>0,\,\zt\in\Sm_{\vr}$, and
which satisfies the following properties; by compactness of
$\Sm_{\vr}$, there are a sufficiently large finite number of
points $\zt_1,\zt_2,\cdots,\zt_{N_0}\in\Sm_{\vr}$ such that for
$\el=1,2,\cdots,N_0$,

(i) $\sum_{\el=1}^{N_0}\Pi_{\el}(\zt)\equiv 1$ for all
$\zt\in\Sm_{\vr}$,

(ii) $\Xi_{\el}(\zt)=1$ for all $\zt\in \Sm_{\vr}\cap
 B(\zt_{\el};2^{-M/2})$,

(iii) $\Xi_{\el}$ is supported in $\Sm_{\vr}\cap
 B(\zt_{\el};2^{1-M/2})$,

(iv) $\lt|\cD^{\ap}\Pi_{\el}(\xi)\rt|\ls 2^{|\ap|M/2}$ for any
multi-index $\ap$, if $1/2\le\vr(\xi)\le 2$,

(v) $N_0\ls 2^{(n-1)M/2}$ for some sufficiently large fixed $M$ (
to be chosen later ),

\noindent where $B(\zt_0;s)$ denotes the ball in $\BR^d$ with
center $\zt_0\in\Sm_{\vr}$ and radius $s>0$. Then we split the
kernel $\cK^{\dt(p)}_l(x,t)$ into finitely many pieces as follows;
$$ \split &\cK^{\dt(p)}_l(x,t) \\ &=2^{(d+1)l}\sum_{\el=1}^{N_0}\iint
e^{i 2^l|x|\la \f{x}{|x|},\xi\ra\tau+i 2^l t\tau}
(1-\vr(\xi))_+^{\dt(p)}\vp(\xi)\,\Pi_{\el}(\xi)\,\psi(\tau)\,\tau^d\,d\xi\,d\tau
\\
&\fd\sum_{\el=1}^{N_0}\cK^{\dt(p)}_{l\el}(x,t).\endsplit\tag{3.3}$$
Here we observe that for $l\in\BZ$ and $\el=1,2,\cdots, N_0$,
$$\cK^{\dt(p)}_{l\el}(x,t)=2^{(d+1)l}\cK^{\dt(p)}_{0\el}\,(2^l
x,2^l t). \tag{3.4}$$

\proclaim{Lemma 3.3} Suppose that $\vr\in C^{\iy}(\BR^d\s\{0\})$
is a non-radial homogeneous distance function satisfying
$\vr(t\xi)=t\vr(\xi)$. If $\dt(p)=d(1/p-1/2)-1/2$ for $0<p<1$,
then given $n,N\in\BN$, there exists a constant
$C=C(p,n,N,\Sm_{\vr})$ such that for any $(x,t)\in\Gm_{\gm}$,
$$ \split &|\cK_{0\el}^{\dt(p)}(x,t)|+\sum_{|\ap|=n+1}\f{1}{\ap !}\,
|[\cD^{\ap}\cK_{0\el}^{\dt(p)}](x,t)|\\ &\le
\f{C}{(1+|x|)^{\dt(p)+1+\f{d-1}{2}}}
\,\,\varPhi\lt(\f{x}{|x|}\rt)\f{1}{(1+|t|-\gm |x|)^{N}}
\endsplit
$$ where $\ap\in (\BN\cup\{0\})^d$ is a
multi-index satisfying $|\ap|=n+1$.
\endproclaim

\noindent{\it Proof.} We observe that the map
$$\BR_+\times\Sm_{\vr}\to\BR^d_0,\,(\vr,\zt)\mapsto\vr\,\zt=\xi,\,\zt\in\Sm_{\vr}$$
defines polar coordinates with respect to $\vr$ by way of
$$d\xi=\vr^{d-1}\la \zt,\fn(\zt)\ra\,d\vr\,d\sm(\zt),$$ where $d\sm(\zt)$
denotes the surface area measure on $\Sm_{\vr}$ and $\fn(\zt)$ is
the outer unit normal vector to $\Sm_{\vr}$ at $\zt\in\Sm_{\vr}$.
Now fix $\zt_{\el}\in\Sm_{\vr}$. Then the unit sphere $\Sm_{\vr}$
can be parametrized near $\zt_{\el}\in\Sm_{\vr}$ by a map
$$w\mapsto \cP (w),\,w\in \BR^{d-1},\,|w|<1$$ such that $\cP (0)=\zt_{\el}$.
Then there is a neighborhood $\cU_0$ of $\zt_{\el}$ with compact
closure, and a neighborhood $\cV_0$ of the origin in $\BR^{d-1}$
so that the map
$$\cQ :(1/2,3/2)\times \cV_0\to \cU_0,(\vr,w)\mapsto \cQ (\vr,w)=\vr\,\cP (w)\tag{3.5}$$
is a diffeomorphism with $\cQ (1,0)=\zt_{\el}$. The Jacobian of
$\cQ$ is given by
$$\fJ (\vr,w)={\vr}^{d-1}\la \cP (w),\fn(\cP (w)) \ra\,\cR (w),$$
where $\cR (w)$ is positive and
$$[\cR (w)]^2=\det \lt({\lt[\f{d\cP}{dw}\rt]}^* \lt[\f{d\cP}{dw}\rt]\rt).$$
By (3.3), we have that $$ \cK_{0\el}^{\dt(p)}(x,t)=\iint e^{i
|x|\la \f{x}{|x|},\xi\ra\tau+i t\tau}
(1-\vr(\xi))_+^{\dt(p)}\vp(\xi)\,\Pi_{\el}(\xi)\,
\tau^d\,\psi(\tau)\,d\xi\,d\tau. \tag{3.6}$$ Let
$\cI_{\el}^0(x,\tau)$ denote the integral with respect to
$\xi$-variable in (3.6); that is to say,
$$\cI_{\el}^0(x,\tau)=\int_{\BR^d} e^{i |x|\la
\f{x}{|x|},\xi\ra\tau}
(1-\vr(\xi))_+^{\dt(p)}\vp(\xi)\,\Pi_{\el}(\xi)\,d\xi.$$ Applying
generalized polar coordinates that we introduced in the above, we
have that
$$ \split &\cI_{\el}^0(x,\tau) \\
&=\iint  e^{i |x|\la \f{x}{|x|},\vr\,\zt\ra\tau}
(1-\vr)_+^{\dt(p)}\,\vp(\vr\,\zt)\,\Pi_{\el}(\vr\,\zt)
\,\vr^{d-1}\la\zt,\fn(\zt)\ra\,d\vr\,d\sm(\zt)
\\&=\int\lt[\int  e^{i |x|\la \f{x}{|x|},\vr\,\cP(w)\ra\tau}
\,\vp[\cQ(\vr,w)]\,\Pi_{\el}[\cQ(\vr,w)]
\,\fJ(\vr,w)\,dw\rt](1-\vr)_+^{\dt(p)}\,d\vr \\
&\fd\int\cJ^0_{\el}(x,\tau,\vr)\,(1-\vr)_+^{\dt(p)}\,d\vr.
\endsplit \tag{3.7}$$ We note that if $\la\th,\fn(\zt_{\el})\ra<1$,
then we have that
$$\ds\lt.\n_w\la\th,\cP(w)\ra\,\rt|_{w=0}\neq 0.$$ Combining this
with the homogeneity condition on the distance function $\vr$ and
choosing a sufficiently large $M>0$ in (i)$\,\,\sim\,\,$(v), we
may select $\vep_0>0$, a neighborhood $\cU_1$ of $\zt_{\el}$ with
$\supp(\vp\,\Pi_{\el})\subset\overline {\cU_1}\subset\cU_0$, and a
neighborhood $\cV_1$ of the origin in $\BR^{d-1}$ so that (3.5)
satisfies, and such that for all $(w,\vr)\in\overline{\cV_1}\times
[1/2,1]$,
$$|\n_w\la\th,\vr\,\cP(w)\ra|\ge c_0\,\,\,\text{ if
$|\la\th,\fn(\zt_{\el})\ra|\le 1-\vep_0$ }$$ and $$c_1\le
\lt|\f{\pa}{\pa\vr}\la\th,\vr\,\cP(w)\ra\rt|\le c_2\,\,\,\text{ if
$|\la\th,\fn(\zt_{\el})\ra|\ge 1-\vep_0$ }\tag{3.8}$$ for some
$c_0>0$, $c_1>0$, and $c_2>0$. We choose some $\fe\in S^{d-2}$ so
that for all $(w,\vr)\in\overline{\cV_1}\times [1/2,1]$,
$$|\cD^w_{\fe}\la\th,\vr\,\cP(w)\ra|\ge \f{1}{2}\,|\n_w\la\th,\vr\,\cP(w)\ra|\ge
\f{1}{2}\,c_0\,\,\,\text{ if $|\la\th,\fn(\zt_{\el})\ra|\le
1-\vep_0$. }\tag{3.9}$$

If $|\la\th,\fn(\zt_{\el})\ra|\le 1-\vep_0$, we apply the
integration of $\cJ^0_{\el}(x,\tau,\vr)$ by parts with respect to
$w$-variable $N_1$-times to obtain that
$$ \cJ^0_{\el}(x,\tau,\vr)=\int  e^{i |x|\la \f{x}{|x|},\vr\,\cP(w)\ra\tau}
\,\f{(\cD^w_{\fe})^{N_1}\lt(\vp[\cQ(\vr,w)]\,\Pi_{\el}[\cQ(\vr,w)]
\,\fJ(\vr,w)\rt)}{(i\,
|x|\,\cD^w_{\fe}\la\th,\vr\,\cP(w)\ra\,\tau)^{N_1}}\,dw. $$ Then
we integrate the kernel $\cK_{0\el}^{\dt(p)}(x,t)$ by parts with
respect to $\tau$-variable $N$-times to get that
$$ \split \cK_{0\el}^{\dt(p)}(x,t)&=\iiint e^{i |x|\la
\f{x}{|x|},\vr\,\cP(w)\ra\tau+i t\tau}
\f{(\f{\pa}{\pa\tau})^N\{\tau^{d-N_1}\,\psi(\tau)\}}
{[i(t+|x|\la\f{x}{|x|},\vr\,\cP(w)\ra)]^N}\,d\tau \\
&\times
\f{(\cD^w_{\fe})^{N_1}\lt(\vp[\cQ(\vr,w)]\,\Pi_{\el}[\cQ(\vr,w)]
\,\fJ(\vr,w)\rt)}{(i\,|x|\,\cD^w_{\fe}\la\th,\vr\,\cP(w)\ra)^{N_1}}\,dw
\,(1-\vr)_+^{\dt(p)}\,d\vr.
\endsplit\tag{3.10}$$ Thus, by Lemma 3.1, (3.9), and (3.10), we have that for
any $N,N_1\in\BN$, $$|\cK^{\dt(p)}_{0\el}(x,t)|\le
\f{C_{N_1}}{(1+|x|)^{N_1}} \,\,\f{1}{(1+|t|-\gm
|x|)^N}.\tag{3.11}$$

If $|\la\th,\fn(\zt_{\el})\ra|\ge 1-\vep_0$, then we apply the
integration of $\cK_{0\el}^{\dt(p)}(x,t)$ ( in (3.6) and (3.7) )
by parts with respect to $\tau$-variable $N$-times to obtain that
$$ \split \cK_{0\el}^{\dt(p)}(x,t)
&=\iiint e^{i |x|\la\f{x}{|x|},\vr\,\cP(w)\ra\tau+i
t\tau}\f{(\f{\pa}{\pa\tau})^N\{\tau^d\,\psi(\tau)\}}
{[i(t+|x|\la\f{x}{|x|},\vr\,\cP(w)\ra)]^N}\,d\tau \\
&\,\,\times \vp[\cQ(\vr,w)]\,\Pi_{\el}[\cQ(\vr,w)]
\,\fJ(\vr,w)\,dw
\,(1-\vr)_+^{\dt(p)}\,d\vr \\
&=\iiint e^{i |x|\la
\f{x}{|x|},\vr\,\cP(w)\ra\tau}\,(1-\vr)_+^{\dt(p)} \\
&\,\,\times \f{\vp[\cQ(\vr,w)]\,\Pi_{\el}[\cQ(\vr,w)]
\fJ(\vr,w)}{[i(t+|x|\la\f{x}{|x|},\vr\,\cP(w)\ra)]^N}\,d\vr\,
dw\,\,\,e^{i
t\tau}\lt(\f{\pa}{\pa\tau}\rt)^N\{\tau^d\psi(\tau)\}\,d\tau.
\endsplit\tag{3.12}$$
Then it follows from the asymptotic result (see [4]) of (3.12)
with respect to $\vr$-variable that for any $N_1\in\BN$,
$$\cK^{\dt(p)}_{0\el}(x,t)=\sum_{j=0}^{N_1 -1}\int e^{i t\tau}\,\fh_j(x,\tau)
\lt(\f{\pa}{\pa\tau}\rt)^N\{\tau^d\psi(\tau)\}\,d\tau+\cO(|x|^{-N_1}),\tag{3.13}$$
where $$\fh_j(x,\tau)=|x|^{-\dt(p)-1-j}\int_{\Sm_{\vr}} e^{i
|x|\la\f{x}{|x|},\zt\ra\tau}\,\fk_j(\zt,x/|x|)\,d\sm(\zt)$$ and
$\fk_j\in C_0^{\iy}(\cP(\cV_0)\times S^{d-1})$ for
$j=0,1,2,\cdots,N_1 -1$. In particular,
$$\fk_0(\zt,x/|x|)=\Gm(\dt(p)+1)\,e^{-i\pi(\dt(p)+1)/2}\,\f{\Pi_{\el}(\zt)}
{[i(t+|x|\la\f{x}{|x|},\zt\ra)]^N}\,\la\zt,\fn(\zt)\ra\,[\la
x/|x|,\zt\ra]^{-\dt(p)-1}.$$ If we restrict to
$(x,t)\in\Gm_{\gm}$, then the required decay estimate of the
kernel follows from (2.1), (2.2), Lemma 3.1, (3.8), (3.11), and
(3.13). Finally, we can complete the proof by noting the fact that
$\sum_{|\ap|=N}\cD^{\ap}\cK_{0\el}^{\dt(p)}=\sum_{|\ap|=N}(\cD^{\ap}\Psi)
*\cK_{0\el}^{\dt(p)}$ for some Schwartz function
$\Psi\in\fS(\BR^{d+1})$. \qed

\proclaim{Corollary 3.4} Suppose that $\vr\in
C^{\iy}(\BR^d\s\{0\})$ is a non-radial homogeneous distance
function satisfying $\vr(t\xi)=t\vr(\xi)$. If
$\dt(p)=d(1/p-1/2)-1/2$ for $0<p<1$, then given $N\in\BN$, there
exists a constant $C=C(p,N,\Sm_{\vr})$ such that for any
$(x,t)\in\Gm_{\gm}$,
$$ \split &|\cK_0^{\dt(p)}(x,t)|+\sum_{|\ap|=N}\f{1}{\ap !}\,
|[\cD^{\ap}\cK_0^{\dt(p)}](x,t)|\\ &\le
\f{C}{(1+|x|)^{\dt(p)+1+\f{d-1}{2}}}
\,\,\varPhi\lt(\f{x}{|x|}\rt)\f{1}{(1+|t|-\gm |x|)^{N}}
\endsplit
$$ where $\ap\in (\BN\cup\{0\})^d$ is a
multi-index satisfying $|\ap|=N$.
\endproclaim

\subheading{4. The atomic decomposition for $H^p$ spaces and
technical lemmas}

An atom is defined as follws: Let $0<p\le 1$ and $\nu$ be an
integer satisfying $\nu\ge (d+1)(1/p-1)$. A $(p,\nu)$-atom is a
function $\fa$ which is supported on a cube $Q$ with center
$x_0\in\BR^{d+1}$ and which satisfies

(i) $|\fa(x)|\le |Q|^{-1/p}\,\,\,\,$ and $\,\,\,\,$(ii)
$\ds\int_{\BR^{d+1}}\fa(x)\,x^{\ap}\,dx=0,$

\noindent where $\ap=(\ap_1,\ap_2,\cdots,\ap_{d+1})\in
(\BN\cup\{0\})^{d+1}$ is a multi-index satisfying $|\ap|
\fd\ap_1+\ap_2+\cdots +\ap_{d+1} \le\nu$ and $x^{\ap}=x_1^{\ap_1}
x_2^{\ap_2} \cdots x_{d+1}^{\ap_{d+1}}$. If $f=\sum_{i=1}^{\iy}
c_i\fa_i$ where each $\fa_i$ is a $(p,\nu)$-atom and
$\{c_i\}\in\el^p(\BR)$, then $f\in H^p$ and $\|f\|_{H^p}\ls\sum_i
|c_i|^p$, and the converse inequality also holds (see [9]).

First of all, we recall a useful lemma [10] due to Stein,
Taibleson and Weiss on summing up weak type functions in case of
$0<p<1$.

\proclaim{Lemma 4.1} Let $0<p<1$. Suppose that $\{\fh_k\}$ is a
sequence of nonnegative measurable functions defined on a subset
$\Om$ of $\BR^d$ such that
$$\lt|\{x\in\Om :\, \fh_k(x)>\ld\}\rt|\le
\f{A}{\ld^p},\,\,\ld>0,$$ where $A>0$ is a constant. If $\{a_k\}$
is a sequence of positive numbers with $\|\{a_k\}\|_{\el^p}<\iy$,
then we have that
$$\lt|\{x\in\Om :\, \sum_k a_k \fh_k(x)>\ld\}\rt|\le
\f{2-p}{1-p}\,\,\|\{a_k\}\|^p_{\el^p}\, \f{A}{\ld^p},\,\,\ld>0.$$
\endproclaim

Theorem 1.1 is obtained by applying a natural variant of Lemma 4.1
as in the following lemma.

\proclaim{Lemma 4.2} Let $a>0$ and $0<p<1$. Suppose that
$\{\fg_l\}$ is a sequence of measurable functions defined on a
subset $\Om$ of $\BR^d$ such that
$$\lt|\{x\in\Om :\,|\fg_l(x)|>\ld\}\rt|\ls 2^{-a l p}\,\ld^{-p}$$ for $l\in\BN$
and all $\ld>0$. Then we have that $$\lt|\{x\in\Om
:\,\sum_{l\in\BN} |\fg_l(x)|>\ld\}\rt|\ls \ld^{-p}.$$
\endproclaim

For $j\in\BN$ and $l\in\BZ$, we set $$ \split
A_l&=\{(x,t)\in\BR^d\times\BR : 2^l \gm |x|\le 4, 2^l |t|>2\}, \\
B_l&=\{(x,t)\in\BR^d\times\BR : 2^l\gm |x|> 4,\, 2^l\,
|t|\le 2,\, | |t|-\gm |x| |>2^{-l}\},  \\
C_l&=\{(x,t)\in\BR^d\times\BR : 2^l\gm |x|> 4, \,2^l\,
|t|> 2,\, | |t|-\gm |x| |\le 2^{-l}\}, \\
D_l&=\{(x,t)\in\BR^d\times\BR : 2^l\gm
|x|> 4, \,2^l\, |t|> 2,\, | |t|-\gm |x| |> 2^{-l}\} \\
&\bigcap\lt(\{(x,t)\in\BR^d\times\BR : |t|\le 2^{-1}\gm\,
|x|\}\cup\{(x,t)\in\BR^d\times\BR : |t|> 2\gm |x|\}\rt),
\\ E_{jl}&=\{(x,t)\in\BR^d\times\BR : 2^l\gm
|x|> 4, \,2^l\, |t|> 2,\, 2^{-l} 2^{j-1}<| |t|-\gm |x| |\le 2^{-l} 2^j\} \\
&\bigcap\{(x,t)\in\BR^d\times\BR : 2^{-1}\gm\, |x|<|t|\le 2\gm
|x|\}.
\endsplit \tag{4.1}$$ For $j\in\BN$, $l\in\BZ$, $(x,t)\in\BR^d\times\BR$, and $a,b,c\in\BR_+$, we set
$$ \split \cA^c_l(x,t)&=2^{l c}\,
|t|^{-\dt(p)-1}\,\chi_{A_l}(x,t), \\
\cB^b_l(x,t)&= 2^{l b}\,
|x|^{-d/p-N}\,\,\varPhi\lt(\f{x}{|x|}\rt)\,\chi_{B_l}(x,t), \\
\cC^a_l(x,t)&= 2^{l a}\,
|x|^{-d/p}\,\,\varPhi\lt(\f{x}{|x|}\rt)\,\chi_{C_l}(x,t), \\
\cD^b_l(x,t)&= 2^{l b}\,
|x|^{-d/p}\,\,\varPhi\lt(\f{x}{|x|}\rt)\,|t|^{-N}\,\chi_{D_l}(x,t), \\
\cE^a_{jl}(x,t)&= 2^{l a} 2^{-jN}\,
|x|^{-d/p}\,\,\varPhi\lt(\f{x}{|x|}\rt)\,\chi_{E_{jl}}(x,t).
\endsplit \tag{4.2}$$ Then by simple computation we obtain the following
lemma.

\proclaim{Lemma 4.3} Let $0<p<1$ and $\dt(p)=d(1/p-1/2)-1/2$. Then
we have the following estimates;

$(a)$ $\lt|\{(x,t)\in\BR^d\times\BR : \cA^c_l(x,t)>\ld\}\rt|\ls
2^{l h}\,\ld^{-p},\,\,\ld>0,$

$(b)$ $\lt|\{(x,t)\in\BR^d\times\BR : \cB^b_l(x,t)>\ld\}\rt|\ls
2^{l h}\,\ld^{-p},\,\,\ld>0,$

$(c)$ $\lt|\{(x,t)\in\BR^d\times\BR : \cC^a_l(x,t)>\ld\}\rt|\ls
2^{l h}\,\ld^{-p},\,\,\ld>0,$

$(d)$ $\lt|\{(x,t)\in\BR^d\times\BR : \cD^b_l(x,t)>\ld\}\rt|\ls
2^{l h}\,\ld^{-p},\,\,\ld>0,$

$(e)$ $\lt|\{(x,t)\in\BR^d\times\BR : \cE^a_{jl}(x,t)>\ld\}\rt|\ls
2^{-j(Np-1)}\, 2^{l h}\,\ld^{-p},\,\,\ld>0,$ where

$(i)$ $h=(d+1)(p-1)$ for $a=d+1-d/p$, $b=d+1-d/p-N$, and
$c=d-\dt(p)$,

$(ii)$ $h=(d+1+N)p-(d+1)$ for $a=d+1+N-d/p$, $b=d+1-d/p$, and
$c=d+N-\dt(p)$.

\noindent Moreover, if $N>\max\{(d+1)(1/p-1),1/p\}$, then it
easily follows from Lemma 4.1 that $$|\{(x,t)\in\BR^d\times\BR :
\cE^a_l(x,t)>\ld\}|\ls 2^{l h}\,\ld^{-p},\,\,\ld>0,$$ where
$\cE^a_l(x,t)= \sum_{j\in\BN}\cE^a_{jl}(x,t).$
\endproclaim

\noindent{\it Proof.} It easily follows from polar coordinates,
Corollary 2.1, Fubini's theorem, and Chebyshev's inequality. \qed

\subheading{5. Weak type estimates on $H^p(\BR^{d+1})$, $0<p<1$ }

In this section, first of all we shall obtain the uniform weak
type estimates of $\cT^{\dt(p)}\fa\,\,$ on the closed unbounded
conical sector in $\BR^{d+1}$ when $\fa$ is a $(p,N)$-atom with
$N\ge (d+1)(1/p-1)$. Then we complete its weak type estimates on
the whole space $\BR^{d+1}$ by applying the translation invariance
of the operator $\cT^{\dt}$ and $H^p$ spaces.

\proclaim{Proposition 5.1} Suppose that $\vr\in
C^{\iy}(\BR^d\s\{0\})$ is a non-radial homogeneous distance
function satisfying $\vr(t\xi)=t\vr(\xi)$ whose unit sphere
$\Sm_{\vr}$ is a convex hypersurface of finite type, and let
$0<p<1$ be given. If $\fa$ is a $(p,N)$-atom with $N\ge
(d+1)(1/p-1)$ defined on $\BR^{d+1}$, then there exists a constant
$C=C(d,p)$ such that
$$\lt|\{(x,t)\in\BR^{d+1} :|\cT^{\dt(p)}\fa(x,t)|\,\chi_{\Gm_{\gm}}(x,t)>\ld\}\rt|\le C
{\ld}^{-p},\,\,\,\ld>0.$$
\endproclaim

\noindent{\it Proof.} Since $\cT^{\dt(p)}$ is translation
invariant, we may assume that $\fa$ is supported in a cube $Q$ of
diameter $\frak d>0$ centered at the origin. We observe that
$$ \split &|\{(x,t)\in\BR^{d+1} : |\cT^{\dt(p)}\fa(x,t)|\,\chi_{\Gm_{\gm}}(x,t)> \ld
\}| \\
&\qquad\quad\le |\{(x,t)\in  Q_*\cap\Gm_{\gm}
:|\cT^{\dt(p)}\fa(x,t)|\,\chi_{\Gm_{\gm}}(x,t)>\ld/2 \}|\\
&\qquad\quad+|\{(x,t)\in Q_*^c\cap\Gm_{\gm} :
|\cT^{\dt(p)}\fa(x,t)|\,\chi_{\Gm_{\gm}}(x,t)>\ld/2 \}|
\endsplit \tag{5.1}$$
where $Q_*$ is the cube concentric with $Q$ and with sides of
twice the length, and we will show that each term is bounded by
$C\,\ld^{-p}$.

Suppose $(x,t)\in  Q_*\cap\Gm_{\gm}$. By Plancherel theorem and
H\"older's inequality with $p/2 + 1/q = 1$, we have that
$$\iint_{Q_*\cap\Gm_{\gm}} |\cT^{\dt(p)}\fa(x,t)|^p \,dx\,dt\le C\,||\cT^{\dt(p)}\fa||_{L^2(\BR^{d+1})}^p
\,|Q_*|^{1/q}\le C.$$ Hence, by Chebyshev's inequality, we have
that for all $\ld>0$,
$$|\{(x,t)\in Q_* :|\cT^{\dt(p)}\fa(x,t)|\,\chi_{\Gm_{\gm}}(x,t)>\ld/2 \}|\le
C\,\ld^{-p}.\tag{5.2}$$ Next we want to estimate the following
weak type inequality
$$|\{(x,t)\in  Q_*^c :|\cT^{\dt(p)}\fa(x,t)|\,\chi_{\Gm_{\gm}}(x,t)>\ld/2 \}|\le
C\,\ld^{-p},\,\,\,\,\ld>0.\tag{5.3}$$

We first assume that $\fa$ is supported in the cube $Q^0$ of
diameter $1$ centered at the origin. We consider the case
$(x,t)\in (Q^0_*)^c\cap\Gm_{\gm}$. Fix $l > 0$. Since $\fa$ is
supported in the cube $Q^0$ of diameter $1$, it follows from Lemma
2.3, (3.2), (3.4), Lemma 3.2, Corollary 3.4, (4.1), and (4.2) that
$$ \split &|\cT_l^{\dt(p)}\fa(x,t)|\,\chi_{(Q^0_*)^c\cap\Gm_{\gm}}(x,t)  \\
&\le 2^{(d+1)l}\iint_{Q^0}
|\fa(y,s)|\,|\cK^{\dt(p)}_0(2^l (x-y),2^l (t-s))|\,\chi_{(Q^0_*)^c\cap\Gm_{\gm}}(x,t)\,dy\,ds \\
&\le\iint_{Q^0}\cA_l^{d-\dt(p)}(x-y,t-s)\,dy\,ds+\iint_{Q^0}\cB_l^{d+1-d/p-N}(x-y,t-s)\,dy\,ds
\\
&+\iint_{Q^0}\cC_l^{d+1-d/p}(x-y,t-s)\,dy\,ds+\iint_{Q^0}\cD_l^{d+1-d/p-N}(x-y,t-s)\,dy\,ds
\\ &+\iint_{Q^0}\sum_{j\in\BN}\cE_{jl}^{d+1-d/p}(x-y,t-s)\,dy\,ds \\
&\le\cA_l^{d-\dt(p)}(x,t)+\cB_l^{d+1-d/p-N}(x,t)+\cC_l^{d+1-d/p}(x,t)
+\cD_l^{d+1-d/p-N}(x,t) +\sum_{j\in\BN}\cE_{jl}^{d+1-d/p}(x,t),
\endsplit $$ where $N$ is a positive integer satisfying
$N>\max\{(d+1)(1/p-1),1/p\}.$ Thus, summing up over the indices
$j\in\BN$ by using Lemma 4.2 and Lemma 4.3, we obtain that
$$\lt|\{(x,t)\in (Q^0_*)^c :
|\cT_l^{\dt(p)}\fa(x,t)|\,\chi_{\Gm_{\gm}}(x,t)>\ld/4\}\rt|\ls
2^{(d+1)(p-1)l}\ld^{-p},\,\,\ld>0.$$ Adding up over the indices
$l>0$ by using Lemma 4.2 once again, we easily get that
$$\lt|\{(x,t)\in (Q^0_*)^c : \sum_{l>0}|\cT_l^{\dt(p)}\fa(x,t)|\,\chi_{\Gm_{\gm}}(x,t)
>\ld/4\}\rt|\ls\ld^{-p},\,\,\ld>0.\tag{5.4}$$

We now fix $l\le 0$. Let $N\in\BN$ be an integer satisfying
$N-1\le \max\{(d+1)(1/p-1),1/p\}<N$. Then we see that
$(d+1+N)p-(d+1)>0$.  If $(x,t)\in (Q^0_*)^c\cap\Gm_{\gm}$, let
$\cP_{l,x,t}(y,s)$ denote the $(N-1)$-th order Taylor polynomial
of the function $(y,s)\mapsto \cK_l(x-y,t-s)$ expanded near the
origin $(0,0)\in\BR^d\times\BR$. Then we have that
$\cP_{l,x,t}(x,t)=2^{(d+1)l}\,\cP_{0,x,t}(2^l x,2^l t)$ for fixed
$l\le 0$. Then it follows from the vanishing moment conditions on
$\fa$, Lemma 2.3, (3.2), (3.4), Lemma 3.2, Corollary 3.4, (4.1),
and (4.2) that
$$ \split
&|\cT_l^{\dt(p)}\fa(x,t)|\,\chi_{(Q_*^0)^c\cap\Gm_{\gm}}(x,t)\\
&=2^{(d+1)l}\iint_{Q^0} |\fa(y,s)|\,|\cK^{\dt(p)}_0(2^l
(x-y),2^l (t-s))- \cP_{0,x,t}(2^l y,2^l s)|\,\chi_{(Q_*^0)^c\cap\Gm_{\gm}}(x,t)\,dy\,ds \\
&\ls 2^{(d+1)l}\int_0^1\iint_{Q^0}\,\sum_{|\ap|=N} \f{1}{\ap
!}\,|[\cD^{\ap}\cK^{\dt(p)}_0](2^l (x-\tau y),2^l (t-\tau
s))|\,\,|2^l
(y,s)|^N\,\chi_{(Q_*^0)^c\cap\Gm_{\gm}}(x,t)\,dy\,ds\,d\tau \\
&\le \int_0^1\iint_{Q^0}\cA_l^{d+N-\dt(p)}(x-\tau y ,t-\tau s
)\,dy\,ds\,d\tau+\int_0^1\iint_{Q^0}\cB_l^{d+1-d/p}(x-\tau y
,t-\tau s)\,dy\,ds\,d\tau  \\
&+\int_0^1\iint_{Q^0}\cC_l^{d+1+N-d/p}(x-\tau y,t-\tau
s)\,dy\,ds\,d\tau +\int_0^1\iint_{Q^0}\cD_l^{d+1-d/p}(x-\tau
y,t-\tau s)\,dy\,ds\,d\tau \\
&+\int_0^1\iint_{Q^0}\sum_{j\in\BN}\cE_{jl}^{d+1+N-d/p}(x-\tau y,t-\tau s)\,dy\,ds\,d\tau \\
&\le\cA_l^{d+N-\dt(p)}(x,t)+\cB_l^{d+1-d/p}(x,t)+\cC_l^{d+1+N-d/p}(x,t)
+\cD_l^{d+1-d/p}(x,t) +\sum_{j\in\BN}\cE_{jl}^{d+1+N-d/p}(x,t).
\endsplit $$
Therefore summing up over the indices $j\in\BN$ by using Lemma 4.2
and Lemma 4.3 leads us to obtain the following weak type estimate
$$\lt|\{(x,t)\in (Q^0_*)^c :
|\cT_l^{\dt(p)}\fa(x,t)|\,\chi_{\Gm_{\gm}}(x,t)>\ld/4\}\rt|\ls
2^{[(d+1+N)p-(d+1)]l}\ld^{-p},\,\,\ld>0.$$ Adding up over the
indices $l\le 0$ by using Lemma 4.2 once again, we easily obtain
that
$$\lt|\{(x,t)\in (Q^0_*)^c : \sum_{l\le 0}|\cT_l^{\dt(p)}\fa(x,t)|\,\chi_{\Gm_{\gm}}(x,t)
>\ld/4\}\rt|\ls\ld^{-p},\,\,\ld>0.\tag{5.5}$$ Suppose now that $\fa$ is an
arbitrary $(p,N)$-atom ( $N\ge (d+1)(1/p-1)$ ) supported in a cube
$Q$ of diameter $\frak d>0$ centered at $(x_0,t_0)\in\BR^{d+1}$.
Let $\fb(x,t) = {\frak d} ^{(d+1)/p}\,\fa(\frak d\,(x-x_0), \frak
d\,(t-t_0))$. Since $\cT^{\dt(p)}$ is translation invariant,
without loss of generality we may assume that $(x_0,t_0)=(0,0)$.
Then $\fb$ is an atom supported in the cube $Q^0$ of diameter $1$
centered at the origin $(0,0)\in\BR^d\times\BR$ . This implies
that
$$ \split \cT_l^{\dt(p)}\fa(x,t)&= {\frak d}^{-(d+1)/p} \int_{\BR}\int_{\BR^d}
      \fb\lt(\f{x-y}{\frak d},\f{t-s}{\frak d}\rt)\,\cK^{\dt(p)}_l(y,s)\,dy\,ds \\
&= {\frak d}^{-(d+1)(1/p-1)}\,\fb * \cK^{\dt(p)}_l(\frak d\cdot,
\frak d \cdot)\lt (\f{x}{\frak d},\f{t}{\frak d}\rt).\endsplit
\tag{5.6} $$ Repeating the same arguments used in (5.4) and (5.5)
in terms of (5.6), we obtain the weak type estimate (5.3) given in
the above. Hence by (5.1), (5.2), and (5.3) we complete the proof.
\qed

{\it Proof of Theorem 1.1.} Let $f=\sum_{i=1}^{\iy}c_i \fa_i\in
H^p(\BR^{d+1})$ where $\fa_i's$ are $(p, N)$-atom $(N \ge
(d+1)(1/p-1))$. Then we see that
$$\|f\|_{H^p}\sim\sum_{i=1}^{\iy}|c_i|^{p} < \infty.$$ By
Proposition 5.1, we obtain that
$$\lt|\{(x,t)\in \Gm_{\gm}  :
|\cT^{\dt(p)}\fa_i(x,t)|>\ld\}\rt|\ls \ld^{-p}, $$ where the
constant $C$ does not depend upon $\ld$ and $\fa_i$. Thus by
applying Stein, Taibleson, and Weiss's lemma ( see [10] ), we have
that
$$\lt|\{(x,t)\in \Gm_{\gm} : |\cT^{\dt(p)}f(x,t)|>\ld\}\rt|\ls\ld^{-p} \, \sum_{i=1}^{\iy}|c_i|^{p}. \tag{5.7}$$

Finally, it remains to show that the inequality (5.7) holds on
$\BR^{d+1}$. For any $R>0$ there must be a ball
$B_{R}(x_{0},t_{0})$ of radius $R$ centered at a point $(x_{0},
t_{0})$ contained in the conical sector $\Gm_{\gm}$. Then it
follows that
$$ \lt|\{(x,t)\in B_{R}(x_{0}, t_{0})  :
|\cT^{\dt(p)}f(x,t)|>\ld\}\rt|  \le \lt|\{(x,t)\in \Gm_{\gm} :
|\cT^{\dt(p)}f(x,t)|>\ld\}\rt|. \tag{5.8}$$ Now we define
$(\tau_{hk}f)(x,y) = f(x-h, y-k)$, $h\in\BR^d, k\in\BR$. Since the
operator $\cT^{\dt(p)}$ is translation invariant and commutes with
the translation operator ( i.e. $\tau_{hk}(\cT^{\dt(p)}f) =
\cT^{\dt(p)}(\tau_{hk}f)$ ), the left-hand side of (5.8) is
rewritten as
$$\split &\lt|\{(x,t)\in B_{R}(x_{0}, t_{0})  :
|\cT^{\dt(p)}f(x,t)|>\ld\}\rt| \\ &\qquad= \lt|\{(x', t')\in
B_{R}(0,0) : |\tau_{(-x_{0})(-y_{0})}
(\cT^{\dt(p)}f)(x', t')|>\ld\}\rt| \\
&\qquad= \lt|\{(x', t')\in B_{R}(0,0)  :
|[\cT^{\dt(p)}(\tau_{(-x_{0})(-y_{0})}f)](x', t')|>\ld\}\rt|.
\endsplit$$
Thus it follows from (5.8) and the fact
$||(\tau_{(-x_{0})(-y_{0})}
 f)||_{H^{p}} = ||f||_{H^{p}}$ that
$$\lt|\{(x, t)\in B_{R}(0,0)  : |  \cT^{\dt(p)}f (x, t)|>\ld\}\rt|
\ls\ld^{-p} \, \sum_{i=1}^{\iy}|c_i|^{p}.\tag{5.9} $$ Therefore,
the inequality (5.9) being uniform in $R > 0$ implies that

$$\lt|\{(x,t)\in \BR^d \times \BR  :
|\cT^{\dt(p)}f(x,t)|>\ld\}\rt|\ls \ld^{-p} \,
\sum_{i=1}^{\iy}|c_i|^{p}\ls\f{\|f\|_{H^p}}{\ld^p}.$$ Hence we
complete the proof. \qed

\noindent{\it Remark.} If $\dt<d(1/p-1/2)-1/2$, it is easily shown
that $\cT^{\dt}$ is not weak type $(p,p)$ on $H^p(\BR^{d+1})$ from
Lemma 4.3. Moreover, the inequalities in Lemma 4.3 are sharp ( see
[10], p.90 ).

\enddocument